\documentclass[12pt]{article}
\usepackage{amssymb}
\parskip=\medskipamount     
\newtheorem{definition}{Definition}[section]

\newtheorem{theorem}[definition]{Theorem}

\newtheorem{remark}[definition]{Remark}
\newtheorem{example}[definition]{Example}

\font\fr=eufm10  scaled \magstep 1   
\font\ddpp=msbm10  scaled \magstep 1  

\def\QED{\hskip0.1em\hfill\null\ \null\nobreak\hfill
\kern3pt\lower1.8pt\vbox{\hrule\hbox   {\vrule\kern1pt\vbox{\kern1.7pt
\hbox{$\scriptstyle   QED$}\kern0.2pt}\kern1pt\vrule}\hrule}}
\def\R{\hbox{\ddpp R}}               
 
\def\PP{\hbox{\ddpp P}}
\def\QQ{\hbox{\ddpp Q}}                              

\def\hfl#1#2{\smash{\mathop{\hbox to 12 mm{\rightarrowfill}}
\limits^{\scriptstyle#1}_{\scriptstyle#2}}}

\begin{document}
\title{MECHANICAL SYSTEMS SUBJECTED TO GENERALIZED NONHOLONOMIC
CONSTRAINTS} 

\author{
Jorge Cort\'es
\and Manuel de Le\'on 
\and David Mart{\'\i}n de Diego
\and Sonia Mart{\'\i}nez 
}

\date{\today}
\maketitle     

Short title: {\em Generalized nonholonomic constraints}\\
1991 MS Classification: 58F05, 70F25, 70F35, 70Hxx

\bigskip

\begin{abstract}
We study mechanical systems subject to constraint functions that can be 
dependent at some points and independent at the rest. Such systems are modelled by means of generalized codistributions. We discuss how the constraint force can transmit an impulse to the motion at the points of dependence and derive an explicit formula to obtain the ``post-impact'' momentum in terms of the ``pre-impact'' momentum. 
\end{abstract}

\newpage

\section{Introduction}

Mechanical systems subjected to nonholonomic constraints
have received a lot of attention in recent years in the literature of
Geometric Mechanics (see
\cite{BatSni,BKMM,FMD,KSB,Koi,Koon1,noholo,Lewis,Marle,Marle1,vdS,VF,Ver} and references therein). Indeed, the dynamics of nonholonomic mechanics
have been described from several approaches: Hamiltonian, Lagrangian
and even Poisson methods have been used.

The constraints which are usually considered in the literature (both linear 
and nonlinear) satisfy a certain regularity condition. That is, they are 
given by a set of independent nonholonomic constraint functions, or, in a global description, by a distribution on the configuration manifold in the
linear case or a submanifold of its tangent bundle in the case of nonlinear
constraints.

However, there is an increasing interest in engineering and robotics in the 
motion of special mechanical systems as, for example, dynamical devices
that locomote with the enviroment via impacts, sudden changes of phase
space, etc. In many cases, the jump of the system's velocity is produced
by an impulse that enforces new constraints on the system. In some cases, 
these systems admit a nice mathematical modelling.

In this paper, we are interested in the following situation. In a local 
description, given a set of constraints $\{\Phi_{1}, \dots, 
\Phi_{m}\}$, one assumes that they become linearly dependent at some points. 
In a global picture, the constraints are given by a generalized 
codistribution with variable rank. One could think of simple examples
that exhibit this kind of behaviour. For instance, imagine a rolling ball
on a surface which is rough on some parts but smooth on the rest. On the
rough parts, it will roll without slipping and, hence, 
nonholonomic linear constraints will be present. However, when the sphere
reaches a smooth part, these constraints will disappear.

The first (and, up to our knowledge, the unique) reference in a
geometrical context for such kind of systems is \cite{Chinos}. In that 
paper, the constraints are provided by a set of global 1-forms on the
configuration manifold and, using the Frobenius theorem, the authors gave a
classification of them according to the existence of 
some special sets that can exert a big influence on the trajectories of 
the system. In particular, the existence of an integral manifold gives
a sort of partial holonomicity with strong implications.
The authors were mainly motivated by problems in motion planning.
However, we are interested, at least in this first approach, in the
geometrical and topological aspects of the problem. In other words,
we are concerned with obtaining the dynamical laws that govern the motion
of the system. 

In consequence, we consider nonholonomic constraints given by a 
generalized codistribution $D$, that is, a codistribution which does not
necessarily have the same rank at all points in the configuration
manifold. This approach leads us to the definition of the concepts of
regular and singular points. It should be noticed that our definitions
are slightly different from those in \cite{Chinos}.
Indeed, the regular points are those 
where the codistribution has locally constant rank. In this sense, the 
generalized codistribution is a regular codistribution, as is commonly 
understood, on the connected components of the set of regular points. The 
singular points are those where 
the codistribution changes its rank. From a dynamical perspective, the
situation on the regular points is already known: we can derive the 
equations of motion following d'Alembert principle and treat them 
making use of the well-developed theory for nonholonomic Lagrangian
systems. 

However, on the singular points the matter is essentially different. The 
classical derivation of the equations of motion no longer works and we 
must solve the problem with other methods. Here, we have adopted a point of 
view strongly inspired by the theory of impulsive mechanics
\cite{Appe,Br,Impul1,Impul2,Impul3,Ka,LT,NF,Pain,R} 
and we use Newton's second law in its integral form \cite{R}. Analyzing
the trajectories which cross the singular set, we have found that, in certain cases, the constraint force can transmit an impulse to the motion. It is precisely the sudden appearence of new constraints (that is, the change of rank of the codistribution) which induces this impulsive character. More precisely, given a motion $q(t)$ crossing the singular set at time $t_0$, we define two vector subspaces of $T^*_{q(t_0)}Q$ as follows: $D^-_{q(t_0)}$ is the limit of all the 1-forms in the codistribution based on $q(t)$, $t<t_0$, when $t \rightarrow t_0^-$. $D^+_{q(t_0)}$ is defined analogously. Our conclusion is that there exists a jump of momentum only if
$D^+_{q(t_0)}$ is not contained in $D^-_{q(t_0)}$ and the ``pre-impact" momentum $p(t_0)_-$ does not satisfy the constraints imposed by $D^+_{q(t_0)}$. In such a case, we propose that the jump is determined by $\Delta p(t_0) \in D^+_{q(t_0)}$ and the condition that the ``post-impact" momentum $p(t_0)_+$ must satisfy the constraints imposed by $D^+_{q(t_0)}$.

To find out about the relation between the theory developed here and the Hamiltonian theory of impact \cite{Vino1,Vino2} applied to this problem, which seems to be a promising possibility, is the object of current research. 

The paper is organized as follows. In Section 2, we introduce the notion
of generalized codistribution, which is just the geometrization of 
constraints with non-constant rank. In Section 3 we review the theory of impulsive forces and impulsive constraints. The ideas exposed here will be helpful in understanding the developments of Section 4, which constitutes the main contribution of this paper, where we study the equations of motion for mechanical systems subjected to generalized constraints and we derive the equations describing the jump of momenta. Finally, in Section 5, some examples are discussed with detail.

\section{Generalized codistributions}

We introduce here the notion of a generalized codistribution. This
notion will be helpful in subsequent sections to model geometrically the 
dynamical systems under consideration, that is, systems subjected to 
constraints which can ``degenerate" at certain points. All
the results in this section are adapted from the ones stated for
generalized distributions in \cite{Va}.

By a {\bf generalized codistribution} we mean a family of linear subspaces 
$D=\{D_q\}$ of the cotangent spaces $T_q^*Q$. Such a codistribution is called 
differentiable if $\forall q \in \hbox{Dom}\, D$, there is a finite number 
of differentiable local 1-forms $\omega_1,...,\omega_l$ defined
on some open neighbourhood $U$ of $q$ such that
$D_{q'}=\hbox{span}\{\omega_1(q'),...,\omega_l(q')\}$ for all $q' \in
U$. 

We define the rank of $D$ at $q$ as $\rho(q)=\hbox{dim} \, D_q$. Given $q_0 \in Q$, if $D$ is 
differentiable, it is clear that $\rho (q)\ge \rho (q_0)$ in a neighbourhood of $q_0$. Therefore, $\rho$ is a lower semicontinuous function. If $\rho$ is a constant function, then $D$ is a codistribution in the usual sense.

For a generalized differentiable codistribution $D$, a point $q \in Q$ will 
be called regular if $q$ is a local
maximum of $\rho$, that is, $\rho$ is constant on an open neighbourhood
of $q$. Otherwise, $q$ will be called a singular point of
$D$. The set $R$ of the regular points of $D$ is obviously open. But,
in addition, it is dense, since if $q_0 \in S = Q \setminus R$,
and $U$ is a neighbourhood of $q_0$, $U$ necessarily contains regular
points of $D$ ($\rho_{|U}$ must have a maximum because it is integer
valued and bounded). Consequently, $q_0 \in \bar{R}$.

Note that in general $R$ will not be connected, as the following simple
example shows:

\begin{example}
{\rm Let us consider $Q = \R^2$ and the general differentiable codistribution
$D_{(x,y)}=\hbox{span}\{ \phi(x)(dx-dy) \}$, where $\phi(x)$ is defined
by 
\[
\phi(x) =
\left\{
\begin{array}{cc}
0 & x \le 0 \\
e^{-\frac{1}{x^2}} & x > 0
\end{array}
\right.
\]
The singular points are those of the $y$-axis, and the connected
components of $R$ are the half-planes $x>0$ (where the rank is $1$) and
$x<0$ (where the rank is $0$).
}
\end{example}

\begin{remark}
{\rm
We note that the notion of singular point defined here is different 
from the one considered in \cite{Chinos}. In that paper, the authors treat the case of generalized constraints given by a globally defined set 
of 1-forms, $\omega_1,...,\omega_l$. Then, they consider the $l$-form
\[
\Omega = \omega_1 \wedge ... \wedge \omega_l \, .
\]
The singular set consists of the points for which $\Omega(q) = 0$, that is, 
the points, $q$, such that $\{ \omega_1(q),...,\omega_l(q) \}$ are linearly 
dependent. Applying this notion to the former example, the set of singular 
points would be the half-plane $\{x \le 0\}$.
}
\end{remark}

Given a generalized codistribution, $D$, we define its annihilator, $D^o$, as the generalized distribution given by
\[
\begin{array}{rrcl}
D^o : & \hbox{Dom} \, D \subset Q & \longrightarrow & TQ \\
& q & \longmapsto & D^o_q = (D_q)^o \, .
\end{array}
\]
Remark that if $D$ is differentiable, $D^o$ is not differentiable, even 
continuous, in general (the corresponding rank function of $D^o$ will
not be lower semicontinuous). In fact, $D^o$ is differentiable if and only if $D$ is a regular codistribution.

We will call $M$ an {\bf integral submanifold} of $D$ if $T_mM$ is
annihilated by $D_m$ at each point $m\in M$. $M$ will be an integral
submanifold of maximal dimension if
\[
T_m M^o = D_m \, , \; \forall m \in M \, .
\]
In particular, this implies that the rank of $D$ is constant along $M$.
A {\bf leaf} $L$ of $D$ is a connected integral submanifold of maximal 
dimension such that every connected integral manifold of maximal
dimension of $D$ which intersects $L$ is an open submanifold of $L$. $D$ will be a {\bf partially integrable} codistribution if for every regular point $q \in R$, there exists one leaf passing through $q$. $D$ will be a {\bf completely integrable} codistribution if there exists a leaf passing through $q$, for every $q\in Q$. In the latter case, the set of leaves defines a general foliation of $Q$. Obviously, any completely integrable codistribution is partially integrable.

$M$ being an integral submanifold of $D$ is exactly the same as being
an integral submanifold of its annihilator $D^o$, and so on.

In Example 2.1, the leaves of $D$ are the half-plane $\{x<0 \}$ and the 
half-lines of slope 1 in the half-plane $\{ x>0\}$. Given any singular point, there is no leaf passing through it. Consequently, $D$ is not a completely integrable codistribution, but it is partially integrable.

\section{Impulsive forces}

In this section, we discuss classical mechanical systems with impulsive 
forces \cite{Appe,Ka,NF,Pain,R,Ver}. This field has traditionally been studied by a rich variety of methods (analytical, numerical and experimental), being a meeting place among physicists, mechanical engineers and mathematicians (for an excellent overview on the subject, see \cite{Br}). Recently, such systems have been brought into the context of Geometric Mechanics \cite{Impul1,Impul2,Impul3,LT}. We will give here a brief review of the classical approach. These ideas will be useful in understanding the behaviour of the constraint forces acting on mechanical systems subject to generalized constraints. Both situations are not the same, but have many points in common, as we will see in the following.

Consider a system of $n$ particles in $\R^3$ such that the particle $r$ has 
mass $m_r$. Introducing coordinates $(q^{3r-2},q^{3r-1},q^{3r})$ for the 
particle $r$, we denote by $Q$ the configuration manifold $R^{3n}$ and by $F_r=(F^{3r-2},F^{3r-1},F^{3r})$ the resultant of all forces acting on the $r^{th}$ particle.

The motion of the particle $r$ in an interval $[t,t']$ is determined by the 
system of integral equations
\begin{equation}\label{NL}
m_r(\dot{q}^k(t')-\dot{q}^k(t)) = \int_{t}^{t'} F^k(\tau) d\tau \, ,
\end{equation}

where $3r-2 \le k \le 3r$ and $k$ is an integer. The integrals of the right-hand side are the components of the {\bf impulse} of the force $F_r$. Equation (\ref{NL}) establishes the relation between the impulse and the momentum change, i.e. ``impulse is equal to momentum change". Equation (\ref{NL}) is a generalized writing of Newton's second law, stated in integral form in order to allow us to consider the case of velocities with finite jump 
discontinuities. This is precisely the case of impulsive forces, which 
generate a finite non-zero impulse at some time instants.

If $F$ is impulsive there exists an instant $t_0$ such that
\begin{equation}\label{IF}
\lim_{t \rightarrow t_0} \int_{t_0}^t F(\tau) d\tau = P \not= 0 \, .
\end{equation}

Equation (\ref{IF}) implies that the impulsive force has an infinite 
magnitude at the point $t_0$, but we are assuming that its impulse $P$
is well defined and bounded. The expression $P \cdot \delta (t_{0}) =
\displaystyle{\lim_{t \rightarrow t_{0}}} \, F(t)$ can be
mathematically seen as a Dirac delta function concentrated at $t_0$.

Now, we will derive the equations for impulsive motion following the discussion in \cite{R}. In the sequel, the velocity vector of the $r^{th}$ particle, $(\dot{q}^{3r-2},\dot{q}^{3r-1},\dot{q}^{3r})$, will be denoted by 
$\dot{q}^r$. Then, the system of integral equations (\ref{NL}) can be
written as
\[
m_r(\dot{q}^r(t_0+\epsilon)-\dot{q}^r(t_0-\epsilon)) =
\int_{t_0-\epsilon}^{t_0 + \epsilon} F_r(\tau) 
d\tau \, .
\]
If we multiply this expression by the virtual displacements at the
point $q(t_0)$, we obtain
\[
(p_r(t_0+\epsilon)-p_r(t_0-\epsilon)) \cdot \delta q^r =
\int_{t_0-\epsilon}^{t_0 +\epsilon} F_r(\tau) d\tau 
\cdot \delta q^r \, .
\]
For the entire system, one has 
\begin{equation}\label{ole1}
\sum_{r=1}^n \left\{ p_r(t_0+\epsilon)-p_r(t_0-\epsilon) -
\int_{t_0-\epsilon}^{t_0+\epsilon} F'_{r}(\tau) 
d\tau \right\} \cdot \delta q^r = \sum_{r=1}^n
\int_{t_0-\epsilon}^{t_0+\epsilon} F''_r(\tau) 
d\tau \cdot \delta q^r \, ,
\end{equation}
where $F'_{r}$ and $F''_r$ are, respectively, the resultant of the
{\it given} forces and of the {\it constraint reaction} forces acting on the $r^{th}$ particle at time $\tau$.

Now, take a local chart $(q^A)$, $1 \le A \le 3n$ on a neighbourhood $U$ of $q(t_0)$ and consider the identification $T_qQ \equiv R^{3n}$, which maps each $v_q \in T_qQ$ to $(v^A)$, such that $v_q=v^A \displaystyle{\left( \frac{\partial}{\partial q^A} \right)_q}$, for each $q\in U$. Let us suppose that the constraints are given on $U$ by the 1-forms $\omega_i=\mu_{iA} dq^A$, $1 \le i \le m$. Then, we have that $\mu_{iA}(q(t))=\mu_{iA}(q(t_0)) + O(t-t_0)$ along the trajectory $q(t)$. As the virtual displacements at the point $q(t)$ satisfy by definition
\[
\sum \mu_{iA}(q(t)) (\delta q (t))^A = 0 \, , \; 1 \le i \le m \, ,
\]
we conclude that $\sum \left( \mu_{iA}(q(t_0)) (\delta q(t))^A + O(t-t_0) \right)= 0$. Therefore, we have that
\[
\delta q^r(t) = \delta q^r (t_0) + O(\epsilon) \, , \; t\in [t_0-\epsilon,t_0+\epsilon] \, ,
\]
that is, the virtual displacements at $q(t)$ can be approximated by the virtual displacements at $q(t_0)$. As a consequence, in the right-hand side of (\ref{ole1}) we can write
\[
\int_{t_0-\epsilon}^{t_0+\epsilon} F''_r(\tau) d\tau \cdot \delta q^r =
\int_{t_0-\epsilon}^{t_0+\epsilon} F''_r(\tau) 
\cdot \delta q^r d\tau = \int_{t_0-\epsilon}^{t_0+\epsilon} F''_r(\tau) 
\cdot \delta q^r (\tau) d\tau + O(\epsilon) \, .
\]
The first term after the last equality is the virtual work done by the constraint forces along the trajectory, and this work is zero since we are considering ideal constraints. The second one goes to zero as $\epsilon$ tends to zero.

In the presence of given impulsive forces acting on $m$ particles, say, at 
time $t_0$, we have 
\[
\lim_{t \rightarrow t_0} \int_{t_0}^t F_{r'}(\tau) d\tau = P_{r'} \not=
0 \; ,  \; 1 \le r' \le m \; . 
\]
Then, taking the limit $\epsilon \rightarrow 0$ in (\ref{ole1}), we obtain the equation for impulsive motion \cite{NF,R}
\begin{equation}\label{funeqn1}
\sum_{r=1}^n \left\{ p_r(t_0)_+ -p_r(t_0)_- - P_{r} \right \} \cdot 
\delta q^r = 0 \, .
\end{equation}
An example in which equation (\ref{funeqn1}) can be applied is when we strike with a cue a billiard ball which is initially at rest. In that case we are exerting an impulsive force that puts the ball into motion. But what happens when the ball collides with the edge of the billiard? What we see is that it bounces, i.e. it suffers again a discontinuous jump in its velocity. The constraint imposed by the wall of the billiard exerts an impulsive force on the ball. When the impulsive force is caused by constraints, such constraints are called {\bf impulsive constraints}. There is a number of different situations in which they can appear. In the following, we examine them.

In the presence of linear constraints of type $\Psi=0$, where $\Psi = b_k(q) \dot{q}^k$ (a situation which covers the case of unilateral holonomic constraints, such as the impact against a wall, and more general types of constraints such as instantaneous nonholonomic constraints), the constraint force, $F = F_k\, dq^k$, is given by

\[
F_k = \mu \hspace{-2pt} \cdot \hspace{-2pt} b_k \, ,
\]

where $\mu$ is a Lagrange multiplier. Then the constraint is impulsive if 
and only if

\[
\lim_{t \rightarrow t_0} \int_{t_0}^t \mu \hspace{-2pt} \cdot
\hspace{-2pt} b_k d\tau = P_k \not= 0 \, ,
\]

for some $k$. The impulsive force may be caused by different circumstances: 
the function $b_k$ is discontinuous at $t_0$, the Lagrange multiplier $\mu$ is discontinuous at $t_0$ or both.

The presence of such constraints does not invalidate equation (\ref{funeqn1}). It merely means that the virtual displacements $\delta q^r$ must satisfy certain additional conditions, which are just those imposed by the constraints. So, in the abscence of impulsive external forces and in the presence of impulsive constraints, we would have

\begin{equation}\label{funeqn2}
\sum_{r=1}^n \Delta p_r(t_0) \cdot \delta q^r = 0 \, ,
\end{equation}

where $\Delta p_r(t_0) = p_r(t_0)_+-p_r(t_0)_-$.

\begin{remark}
{\rm In general, equation (\ref{funeqn2}) is not enough to determine the jump of the momentum. One usually needs additional physical hypothesis, related with elasticity, plasticity, etc. to obtain the post-impact momentum. In this respect, there are two classical approaches, the Newtonian approach and the Poisson approach \cite{Br,St}. The Newtonian approach relates the normal component of the rebound velocity to the normal component of the incident velocity by means of an experimentally determined coefficient of restitution $e$, where $0 \le e \le 1$. Poisson approach divides the impact into compression and decompression phases and relates the impulse in the restitution phase to the impulse in the compression phase.
}
\end{remark}

\begin{remark}
{\rm It could happen that impulsive constraints and impulsive forces to be present at the same time. For example, in the collision between a rigid lamina and an immobile plane surface, we must take into account not only the normal component of the contact force, but also the friction force associated to the contact. It is not innocuous the way the friction is entered into the picture. In fact, the Newton and Poisson approaches have been revealed to be physically inconsistent in certain situations. On the one hand, Newton approach can show energy gains \cite{Ke,St}. On the other hand, Poisson's rule is not satisfactory since non-frictional dissipation does not vanish for perfectly elastic impacts \cite{Br,St}. This surprising consequence of the impact laws is only present when the velocity along the impact surface (slip) stops or reverses during collision, due precisely to the friction. Stronge \cite{St,St2} proposed a new energetically consistent hypothesis for rigid body collisions with slip and friction. It should be noticed that the three approaches are equivalent if slip does not stop during collision and in the perfectly inelastic case ($e=0$).

Recently, a new Newton-style model of partly elastic impacts has been proposed \cite{Ste} which, interestingly, always dissipates energy, unlike the classical formulation of the Newtonian approach discussed in \cite{St}.
}
\end{remark}

In the frictionless case, one can prove the following
\begin{theorem}[Carnot's theorem](\cite{Impul3,R})
The energy change due to impulsive constraints is always a loss of energy. 
\end{theorem}

\section{Mechanical systems subjected to generalized constraints}

In this section, we study the equations of motion for mechanical systems subjected to generalized constraints. Let us consider a mechanical system with Lagrangian function $L: TQ 
\rightarrow \R$, $L(v) = \frac{1}{2}g(v,v) - (U\circ \tau_Q)(v)$, where $g$ is a Riemannian metric on $Q$ and $U$ is a function on the
configuration space $Q$ (the potential). Suppose, in addition, that the
system is subjected to a set of constraints given by a generalized 
differentiable codistribution $D$ on $Q$, that is, we assume that 
$\tau_Q(D) = Q$. The motions of the system are forced to take place 
satisfying the constraints imposed by $D$.

We know that the codistribution $D$ induces a decomposition of $Q$ into 
regular and singular points. We write
\[
Q=R \cup S \, .
\]
Let us fix $R_c$, a connected component of $R$. We can consider the 
restriction of the codistribution to $R_c$, $D_c=D_{|R_c}: R_c 
\subset Q \longrightarrow T^*Q$. Obviously, we have that $D_c$ is a regular 
codistribution, that is, it has constant rank.

Then, let us denote by $D_c^o:R_c \longrightarrow TQ$ the annihilator of 
$D_c$. Now, we can consider the dynamical problem with regular Lagrangian $L$, subjected to the regular codistribution $D_c^o$ and apply the 
well-developed theory for nonholonomic Lagrangian systems  \cite{BKMM,Koon1,noholo,Marle}.

Consequently, our problem is solved on each connected component of $R$. The situation changes radically if the motion reaches a singular point. The rank of the constraint codistribution can vary suddenly and the classical derivation of the equations of motion for nonholonomic Lagrangian systems is no longer valid. Let us explore the behaviour of the system when such a thing occurs.

Consider a trajectory of the system, $q(t)$, which reaches a singular point 
at time $t_0$, i.e. $q(t_0) \in S$, such that $q(t_0-\epsilon,t_0) \subset R$ and $q(t_0,t_0+\epsilon) \subset R$ for sufficiently small $\epsilon>0$. The motion along the trajectory
$q(t)$ is governed by the following equation, which is, as in the
impulsive case, an integral writing of Newton's second law, to
consider possible finite jump discontinuities in the velocities (or the
momenta). That is, at each component
\begin{equation}\label{NL2}
p_A(t')-p_A(t) = \int_{t}^{t'} F_A(\tau) d\tau \, ,
\end{equation}
on any interval $t \le t' < \infty$, where $F$ is the resultant of all the 
forces action on the trajectory $q(t)$. In our case, the unique forces acting are the constraint reaction forces.

The nature of the force can become impulsive because of the change of
rank of the codistribution $D$. We summarize the situations that can be
found in Table \ref{figurita}. On entering the singular set, the rank of 
the codistribution $D$ at the singular point $q(t_0)$ can be the same as at 
the preceding points (Case 1) or can be lower (Cases 2 and 3). In these two 
latter situations, the constraints have collapsed at $q(t_0)$ and this 
induces a finite jump in the constraint force. As the magnitude of the force
is not infinite, there is no abrupt change in the momenta. Consequently, 
in all cases, we find no momentum jumps on entering the singular set.

\medskip

\begin{table}[hbt]
\begin{center}
\begin{tabular}{r|c|c|c|}
\cline{2-4}
& \multicolumn{1}{c|}{$q(t_0-\epsilon)$: preceding points}
& \multicolumn{1}{c|}{$q(t_0)$: singular point}
& \multicolumn{1}{c|}{$q(t_0+\epsilon)$: posterior points} \\
\hline
\multicolumn{1}{|r|}{Case 1}
& $\rho=r$ & $\rho_0=r_0=r$ & $\rho>r$ \\
\hline
\multicolumn{1}{|r|}{Case 2}
& $\rho=r$ & $\rho_0=r_0<r$ & $\rho=r_0$ \\
\hline
\multicolumn{1}{|r|}{Case 3}
& $\rho=r$ & $\rho_0=r_0<r$ & $\rho>r_0$ \\
\hline
\end{tabular}
\end{center}
\caption{Possible cases. The rank of $D$ is denoted by $\rho$}\label{figurita}
\end{table}

\medskip

On leaving the singular set, the rank of $D$ at the posterior points can be
the same as at $q(t_0)$ (Case 2) or can be higher (Cases 1 and 3). In
Case 2 nothing special occurs. In Cases 1 and 3, the trajectory must 
satisfy, immediately after the point $q(t_0)$,
{\bf additional constraints} which were not present before. It is in
this sense that we affirm that the constraint force can become
impulsive: if the motion which passes through the singular set and tries
to enter the regular one again does not satisfy the new constraints,
then it experiences a jump of its momentum, due to the presence of the
constraint force. In this way, the new values of the momentum
satisfy the constraints. But one has to be careful: the impulsive
force will act just on leaving $S$, on the regular set. Consequently, we must take into account the virtual displacements associated to the posterior regular points. The underlying idea of the mathematical derivation of the momentum jumps in Section 4.1 is the following: to take an infinitesimal posterior point $q(t)$ to $q(t_0)$, to forget for a moment the presence of the constraints on the path $q(t_0,t)$ and to derive the momentum jump at $q(t)$ due to the appearence of the additional constraints. Afterwards, to make a limit process $t \rightarrow t_0$, cancelling out the interval $(t_0,t)$ where we ``forgot" the constraints. In any case, we will make the convention that the jump happens at $q(t_0)$.

We illustrate the above discussion in the following example.

\begin{example}\label{exa}
{\rm Consider a particle in the plane subjected to the constraints imposed
by the generalized codistribution in Example 2.1. The Lagrangian
function is
\[
\begin{array}{rrcl}
L: & T\R^2 & \longrightarrow & \R \\
& (x,y,\dot{x},\dot{y}) & \longmapsto &
\displaystyle{\frac{1}{2}}m(\dot{x}^2+\dot{y}^2) \, .
\end{array}
\]
On the half-plane $R_1=\{x<0\}$ the codistribution is zero and the
motion is free. Consequently the trajectories are
\begin{eqnarray*}
x & = & \dot{x}_0 t + x_0 \, , \\
y & = & \dot{y}_0 t + y_0 \, .
\end{eqnarray*}

If the particle starts its motion with initial conditions $x_0=-1$,
$y_0=1$, $\dot{x}_0=1$, $\dot{y}_0=0$, after a time $1$, it reaches the
singular set $S = \{x=0\}$. If the motion crosses the $y$-axis, something 
abrupt occurs on entering the half-plane $R_2=\{x>0 \}$, where the 
codistribution is no longer zero and, indeed, imposes the additional 
constraint $\dot{x}=\dot{y}$ (Case 1). We know that the integral manifolds of
$D$ on $R_2$ are half-lines of slope 1, so the particle suffers a
finite jump in the velocity on going through the singular part in order
to adapt its motion to the prescribed direction.

If, on the contrary, the particle starts on $R_2$, say with initial
conditions $x_0=1$, $y_0=1$, $\dot{x}_0=-1$, $\dot{y}_0=-1$, after a
certain time, it reaches the set $S$. On crossing it, nothing special
happens, because the particle finds less contraints to fulfill, indeed,
there are no constraints (Case 2). Its motion on $R_1$ is free, on a straight line of slope 1 and with constant velocity equal to the one at the 
singular point of crossing.
}
\end{example}

\subsection{Momentum jumps}

Now, we derive a formula, strongly inspired by the theory of impulsive
motion, for the momentum jumps which can occur due to the changes
of rank of the codistribution $D$ in Cases 1 and 3.

At $q(t_0)$ we define the following vector subspaces of
$T^*_{q(t_0)} Q$
\begin{eqnarray*}
D^-_{q(t_0)}&=&\{\alpha \in T^*_{q(t_0)} Q\; /\; \exists
\tilde{\alpha}: (t_0-\epsilon, t_0)\rightarrow T^*Q,
\tilde{\alpha}(t)\in D_{q(t)}\  \hbox{and}\ \lim_{t\rightarrow
t_0^-}\tilde{\alpha}(t)=\alpha\}\\ 
D^+_{q(t_0)}&=&\{\alpha \in T^*_{q(t_0)} Q\; /\; \exists
\tilde{\alpha}: (t_0, t_0+\epsilon)\rightarrow T^*Q,
\tilde{\alpha}(t)\in D_{q(t)}\  \hbox{and}\ \lim_{t\rightarrow
t_0^+}\tilde{\alpha}(t)=\alpha\}\\ 
\end{eqnarray*}
From the definition of $D^-_{q(t_0)}$ and $D^+_{q(t_0)}$ we have that
\[
(D^-_{q(t_0)})^{\perp}=\lim_{t\rightarrow t_0^-}(D_{q(t)})^{\perp}\
\hbox{ and }\  
(D^+_{q(t_0)})^{\perp}=\lim_{t\rightarrow t_0^+}(D_{q(t)})^{\perp}
\]
where $ ^{\perp}$ denotes the orthogonal complement with respect to the bilinear form induced by the metric $g$ on the cotangent space $T^*_{q(t_0)}Q$, and the limits $(D^{\perp})^{-}$ and $(D^{\perp})^{+}$ are defined as in the case of $D^{-}$ and $D^{+}$. In the following, we shall not make a notational distinction between the metric $g$ and the induced bilinear form on $T^*_{q(t_0)}Q$. In each case, the precise meaning should be clear from the context.

Since $D$ is a differentiable codistribution then
\[
D_{q(t_0)}\subseteq D^-_{q(t_0)}
\ \hbox{ and }\ D_{q(t_0)}\subseteq D^+_{q(t_0)} \, .
\]

Along the interval $[t_0,t]$, we have
\[
p_A(t)-p_A(t_0) = \int_{t_0}^t F_A(\tau)\, d\tau \; .
\]
Multiplying by the virtual displacements at the point $q(t)$, we obtain
\begin{equation}\label{yoyo}
(p_A(t)-p_A(t_0)) \cdot \delta q^A_{\;|q(t)} = 
\int^{t}_{t_0} F_A(\tau)\, d\tau \cdot \delta 
q^A_{\;|q(t)} \; .
\end{equation}

Summing in $A$, we get
\begin{equation}\label{funeqn}
\sum_{A=1}^n (p_A(t)-p_A(t_0)) \cdot \delta q^A_{\;|q(t)} = 
\sum_{A=1}^n \int_{t_0}^{t} F_A(\tau) d\tau \cdot \delta 
q^A_{\;|q(t)} \, .
\end{equation}

Since we are dealing with ideal constraints, the virtual work vanishes, that 
is 
\[
\sum_{A=1}^n\int^t_{t_0}F_A(\tau)\delta q^A_{\;|q(\tau)}\, d\tau=0 \; .
\]
If $t$ is near $t_0$, then $\tau$ is close to $t$, and $q(\tau)$ is near 
$q(t)$, so $\delta q$ remains both nearly constant and nearly equal to its 
value at time $t$ throughout the time interval $(t_0, t]$, in the same way we exposed in Section 3. Therefore,  
\[
\sum_{A=1}^n\left(\int^t_{t_0}F_A(\tau)
\, d\tau\right)\cdot\delta q^A_{\;|q(t)} = \sum_{A=1}^n\int^t_{t_0}
F_A(\tau)\delta q^A_{\;|q(\tau)}\, d\tau + O(t-t_0) = O(t-t_0) \, . 
\]

Consequently, equation (\ref{funeqn}) becomes
\begin{equation}\label{funeqnw}
\sum_{A=1}^n (p_A(t)-p_A(t_0)) \cdot \delta q^A_{\;|q(t)} = O(t-t_0)
\; .
\end{equation}

Taking limits we obtain
\[
\lim_{t\rightarrow t_0^+} \left(\sum_{A=1}^n(p_A(t)-p_A(t_0))\cdot
\delta q^A_{\;|q(t)} \right) = 0
\]
which implies
\begin{equation}
\sum_{A=1}^n(p_A(t_0)_+-p_A(t_0))\lim_{t\rightarrow t_0^+}
\delta q^A_{|q(t)}=0
\end{equation}
or, in other words,
\begin{equation}\label{ase}
(p_A(t_0)_+-p_A(t_0))\, dq^A \in \lim_{t\rightarrow t^+_0} 
D_{q(t)} = D_{q(t_0)}^+ \, .
\end{equation}

\paragraph{Conclusion:} Following the above discussion, we will 
deduce the existence of jump of momenta depending on the relation between 
$D^-_{q(t_0)}$ and $D^+_{q(t_0)}$. The possible cases are shown in 
Table \ref{pepe}.

\begin{table}[hbt]
\begin{center}
\begin{tabular}{|r|c|}
\hline
$D^+_{q(t_0)} \subseteq D^-_{q(t_0)}$&  there is no jump of momenta\\
\hline
$D^+_{q(t_0)} \not \subseteq D^-_{q(t_0)}$&    possibility of jump of momenta\\
\hline
\end{tabular}
\end{center}
\caption{The two cases}\label{pepe}
\end{table}

In the second case in Table \ref{pepe}, we have a jump of momenta if the ``pre-impact" momentum $p(t_0)_-=p(t_0)$ does not satisfy the constraints imposed by $D^+_{q(t_0)}$, that is
\[
{p_A(t_0)_-dq^A} \notin (D^+_{q(t_0)})^\perp \; .
\]

Our proposal for the equations which determine the jump is then
\[
\left\{
\begin{array}{l}
{(p_A(t_0)_+-p_A(t_0)_-)\, dq^A} \in D^+_{q(t_0)}\\
{p_A(t_0)_+\, dq^A} \in (D^+_{q(t_0)})^{\perp} \, .
\end{array}
\right.
\]

The first equation has been derived above (cf. (\ref{ase})) from the generalized writing of Newton's second law (\ref{NL}). The second equation simply encodes the fact that the ``post-impact" momentum must satisfy the new constraints imposed by $D^+_{q(t_0)}$.

\begin{remark}
{\rm In Cases 1 and 3, the virtual displacements at $q(t_0)$ are
radically different from the ones at the regular posterior points,
because of the change of rank. From a dynamics point of view, these are
the ``main'' ones, since it is on the regular set where an
additional constraint reaction force acts. As we have seen, the momentum jump
happens on just leaving $S$, due to the presence of this additional constraint
force on the regular set. Note that with the procedure we have just
derived, we are taking into account precisely the virtual displacements
at the regular posterior points, and not those of $q(t_0)$. If we took
the virtual displacements at $q(t_0)$ and multiply by them in
(\ref{yoyo}), we would obtain non-consistent jump conditions. This is easy to
see, for instance, in Example \ref{exa}.
}
\end{remark}

An explicit derivation of the momentum jumps for Cases 1 and 3 would
be as follows. Let $m$ be the maximum between $\rho=r$, the rank at the
regular preceding points, and $\rho=s$, the rank at the regular posterior
points. Then there exists a neighbourhood $U$ of $q(t_0)$ and 1-forms 
$\omega_1,...,\omega_m$ such that
\[
D_q = \, \hbox{span} \, \{\omega_1(q),...,\omega_m(q) \} \; , \;
\forall q \in U \, . 
\]
Let us suppose that $\omega_1,...,\omega_s$ are linearly independent at
the regular posterior points (if not, we reorder them). 
Obviously, at $q(t_0)$, these $s$ 1-forms are linearly dependent.
In the following, we will denote by $\omega_i$ the 1-form evaluated at
$q(t)$, ($t$ time immediately posterior to $t_0$) i.e. $\omega_i \equiv
\omega_i(q(t))$, in order to simplify notation.

Since the Lagrangian is of the form $L=T-U$, where $T$ is 
the kinetic energy of the Riemannian metric $g$, that is, 
$L=\frac{1}{2} g_{AB}\dot{q}^A\dot{q}^B-U(q)$, then we have that
\begin{equation}\label{formulita2}
\omega_{jA}(q(t))\dot{q}^{A}(t) = \sum_{A,B} \omega_{jA} 
g^{AB}p_B(t) = 0  \, , \; j=1,...,s \; .
\end{equation}

Using the metric $g$ we have the following decomposition
\[
T^*_qQ = D_q\oplus D_{q}^{\perp},\quad q \in Q \; .
\]
The two complementary projectors associated to this decomposition are 
\begin{eqnarray*}
&&{\cal P}_q: T^*_qQ\longrightarrow  D_q^{\perp}\\
&&{\cal Q}_q: T_q^*Q\longrightarrow  D_q
\end{eqnarray*}
The projector ${\cal P}_q$ is given by 
\[
{\cal P}_q(\alpha_q)=\alpha_q-{\cal C}^{ij}\alpha_q(Z_i)\omega_j,\quad
\alpha_q\in T^*Q 
\]
where 
\[
Z_i = g^{AB}\omega_{iB} \left( \frac{\partial}{\partial q^A} \right)_q
\]
and ${\cal C}^{ij}$ are the entries of the inverse matrix of ${\cal C}$, 
the symmetric matrix with entries ${\cal C}_{ij} = \omega_{iA} g^{AB} 
\omega_{jB}$, or ${\cal C} = \omega g^{-1} \omega^{T}$ with the obvious 
notations.  

By definition
\[
p_A(t_0)_+dq^A|_{q(t_0)} = \lim_{t\rightarrow t^+_0}(p_A(t)dq^A|_{q(t)}) \, .
\]
From (\ref{formulita2}), ${\cal P}_{q(t)}(p_A(t)dq^A|_{q(t)})) = 
p_A(t)dq^A|_{q(t)}$ and then  
\begin{equation}\label{aux}
p_A(t_0)_+dq^A|_{q(t_0)} =
\left(\lim_{t\rightarrow t_0^+}{\cal
P}_{q(t)}\right)(p_A(t_0)_+dq^A|_{q(t_0)}) \in (D^+_{q(t_0)})^\perp \, .
\end{equation}
Combining (\ref{ase}) and (\ref{aux}), we obtain
\[
p_A(t_0)_+dq^A|_{q(t_0)}=\left(\lim_ {t\rightarrow t_0^+}{\cal
P}_{q(t)}\right)\left[ p_A(t_0)_-dq^A|_{q(t_0)}\right]  \; .
\]
In coordinates, this can be expressed as
\begin{equation}\label{formulita3}
p_A(t_0)_+ = p_A(t_0)_- - \lim_{t\rightarrow t_0^+}\left(\sum_{i,j,A,B} 
{\cal C}^{ij}
\omega_{jB} g^{BC}  
\omega_{iA}\right)\Big|_{\;q(t)} p_C(t_0)_- \, , \; A=1,...,n \; .
\end{equation}

Equation (\ref{formulita3}) can be written in matrix form as follows
\begin{equation}\label{formulita4}
p(t_0)_+=\left(Id - \lim_{t \rightarrow t_0^+} (\omega^{T} {\cal C}^{-1} 
\omega g^{-1})_{|q(t)}\right) p(t_0)_- \, .
\end{equation}

With the derived jump rule, we are able to prove the following version of Carnot's theorem for generalized constraints.

\begin{theorem}
The kinetic energy will only decrease by the application of the jump rule (\ref{formulita4}).
\end{theorem}
{\bf Proof:} We have that
\begin{eqnarray*}
g\left(p(t_0)_+,p(t_0)_+\right) &=& g\left( \left( \lim_{t\rightarrow t_0^+} {\cal P}_{q(t)} \right) p(t_0)_- ,p(t_0)_- - \left( \lim_{t\rightarrow t_0^+} {\cal Q}_{q(t)} \right) p(t_0)_-\right) \\
&=& g\left(\left( \lim_{t\rightarrow t_0^+} {\cal P}_{q(t)}\right) p(t_0)_- ,p(t_0)_-\right) \\
&=& g\left(p(t_0)_- ,p(t_0)_-\right) - g\left(\left(\lim_{t\rightarrow t_0^+} {\cal Q}_{q(t)}\right) p(t_0)_- ,p(t_0)_-\right) \, .
\end{eqnarray*}
Since $\displaystyle{g\left( \left( \lim_{t\rightarrow t_0^+} {\cal Q}_{q(t)}\right) p(t_0)_- ,p(t_0)_-\right) = g\left( \left(\lim_{t\rightarrow t_0^+} {\cal Q}_{q(t)}\right) p(t_0)_- ,\left(\lim_{t\rightarrow t_0^+} {\cal Q}_{q(t)}\right) p(t_0)_-\right) \ge 0}$, we can conclude that
\[
\frac{1}{2} g(p(t_0)_+,p(t_0)_+) \le \frac{1}{2} g(p(t_0)_- ,p(t_0)_-) \, .
\]
\QED

In fact, the jump rule (\ref{formulita4}) has the following alternative interpretation. Let $p \in D^+_{q(t_0)}$ and observe that
\begin{eqnarray*}
g(p-p(t_0)_-,p-p(t_0)_-) &=& g(p(t_0)_-,p(t_0)_-) + g(p,p-2 p(t_0)_-) \\
&=& g(p(t_0)_-,p(t_0)_-) + g\left( p,p-2 \left( \lim_{t\rightarrow t_0^+} {\cal P}_{q(t)} \right) p(t_0)_- \right) \, .
\end{eqnarray*}
Now, note that the covector $p=\left( \lim_{t\rightarrow t_0^+} {\cal P}_{q(t)} \right) p(t_0)_- \in D^+_{q(t_0)}$ is such that the expression $g(p-p(t_0)_-,p-p(t_0)_-)$ is minimized among all the covectors belonging to $D^+_{q(t_0)}$. Therefore, the derived jump rule (\ref{formulita4}) can be stated as follows: the ``post-impact" momenta $p(t_0)_+$ is such that the kinetic energy corresponding to the difference of the ``pre-impact" and  ``post-impact" momentum is minimized among all the covectors satisfying the constraints. This is an appropriate version for generalized constraints of the well-known jump rule for perfectly inelastic collisions \cite{Mo}. This is even more clear in the holonomic case, as is shown in Subsection 4.2.

\begin{remark}
{\rm 
So far, we have been dealing with impulsive constraints. More 
generally, we can consider the presence of external impulsive forces 
associated to external inputs or controls. 
Then, equation (\ref{ase}) must be modified as follows 
\begin{equation}\label{ase1}
(p_A(t_0)_+-p_A(t_0)_--P'_A(t_0))\, dq^A_{\; |q(t_0)}\in
\lim_{t\rightarrow t^+_0} D_{q(t)} = D^+_{q(t_0)}
\end{equation}
where $P'_A(t_0)$, $1\leq A\leq n$, are the external impulses at time $t_0$. 
Observe that if $q(t_0)$ is a regular point then  
\begin{equation}\label{formulita31}
p_A(t_0)_+ =p_A(t_0)_-+P'_A(t_0)-\sum_{i,j,A,B} {\cal C}^{ij} \omega_{jB} g^{BC} 
\omega_{iA}P'_C(t_0) \, , \; A=1,...,n \, ,
\end{equation}
and, if $q(t_0)$ is a singular point, we have
\begin{equation}\label{formulita32}
p_A(t_0)_+ = p_A(t_0)_-+P'_A(t_0)- \lim_{t\rightarrow t_0^+} 
\left(\sum_{i,j,A,B} {\cal C}^{ij}
\omega_{jB} g^{BC}  
\omega_{iA}\right)\Big|_{\;q(t)} (p_C(t_0)_-+P'_C(t_0)) \, , \; A=1,...,n \, .
\end{equation}
}
\end{remark}

\subsection{The holonomic case}

We show in this  section a meaningful interpretation of the proposed jump rule (\ref{formulita4}) in case the codistribution $D$ is partially integrable.

Let us consider a trajectory $q(t)\in Q$ which reaches a singular point $q(t_0) \in S$ and falls in either Case 1 or Case 3. Since $Q=\bar{R}$, we have that $q(t_0) \in \bar{L}$, where $L$ is the leaf of $D$ which contains the regular posterior points of the trajectory $q(t)$. On leaving $q(t_0)$, we have seen that the trajectory suffers a finite jump in its momentum in order to satisfy the constraints imposed by $D$, which in this case implies that the trajectory after time $t_0$ belongs to the leaf $L$. Consequently, the jump can be interpreted as a perfectly inelastic collision against the ``wall" represented by the leaf $L$!

Let us see it revisiting Example \ref{exa}.

\begin{example}
{\rm Consider again the situation in Example \ref{exa}. If the motion of the particle starts on the left half-plane going towards the right one, then it is easy to see that $D^-_{(0,y)} = \{0\}$ and $D^+_{(0,y)}=\hbox{span} \{dx-dy\}$. As $D^+_{(0,y)} \not \subseteq D^-_{(0,y)}$, a jump of momenta is possible. In fact, if the ``pre-impact" velocity $(\dot{x}_0,\dot{y}_0)$ does not satisfy $\dot{x}_0 = \dot{y}_0$, the jump occurs and is determined by $\Delta v(t_0) \in D^+_{(0,y)}$ and $\dot{x}(t_0^+)=\dot{y}(t_0^+)$. 
Consequently, we obtain
\begin{eqnarray*}
\dot{x}(t_0^+) &=& \frac{\dot{x}_0+\dot{y}_0}{2} \, , \\
\dot{y}(t_0^+) &=& \frac{\dot{x}_0+\dot{y}_0}{2} \, .
\end{eqnarray*}

We would have obtained the same result if we had considered that our particle hits, in a perfectly inelastic collision, against the ``wall" represented by the half-line of slope 1 contained in $\{x >0 \}$ passing through the point $(0,y)$.   

If the particle starts on the right half-plane towards the left one, the 
roles are reversed and $D^-_{(0,y)}=<dx-dy>$, $D^+_{(0,y)}=\{0\}$. We have 
that $D^+_{(0,y)} \subseteq D^-_{(0,y)}$ and therefore there is no jump.
}
\end{example}

\section{Examples}

Next, we are going to develop two examples illustrating the above
discussion. First, we treat a variation of the classical example of the rolling sphere \cite{NF,R}. Secondly, we take one example from Chen, Wang, Chu and Chou \cite{Chinos}.

\subsection{The rolling sphere} 

Consider a homogeneous sphere rolling on a plane.
The configuration space is $Q=\R^2 \times SO(3)$: $(x,y)$ denotes the
position of the center of the sphere and $(\varphi, \theta, \psi)$
denote the Eulerian angles.  

Let us suppose that the plane is smooth if $x<0$ and absolutely
rough if $x>0$. On the smooth part, we assume that 
the motion of the
ball is free, that is, the sphere can slip. But if it reaches the rough
half-plane, the sphere begins rolling without slipping, because of the
presence of the constraints imposed by the roughness. We are interested
in  knowing the trajectories of the sphere and, in particular, the
possible  changes in its dynamics because of the crossing from one
half-plane to the other.

The kinetic energy of the sphere is
\[
T = \frac{1}{2}
\left(\dot{x}^2+\dot{y}^2+k^2(\omega_x^2+\omega_y^2+
\omega_z^2)\right) \, ,
\]
where $\omega_{x}$, $\omega_{y}$ and $\omega_{z}$ are the angular
velocities with respect to the inertial frame, given by
\begin{eqnarray*}
\omega_{x} & = & \dot{\theta} \cos \psi + \dot{\varphi} \sin \theta
\sin \psi \, , \\
\omega_{y} & = & \dot{\theta} \sin \psi - \dot{\varphi} \sin \theta
\cos \psi \, , \\
\omega_{z} & = & \dot{\varphi} \cos \theta + \dot{\psi} \, .
\end{eqnarray*}
The potential energy is not considered here since it is constant. 

The condition of rolling without sliding of the sphere when $x>0$
implies that the point of contact of the sphere and the plane has zero
velocity  
\begin{eqnarray*}
\phi^1&=&\dot{x}-r\omega_y=0 \, , \\
\phi^2&=&\dot{y}+r\omega_x=0 \, ,
\end{eqnarray*}
where $r$ is the radius of the sphere.

Following the classical procedure, we introduce quasi-coordinates
``$q_{1}$'', ``$q_{2}$'' and ``$q_{3}$'' such that ``$\dot{q}_{1}$''$ =
\omega_x$, ``$\dot{q}_{2}$''$= \omega_y$ and ``$\dot{q}_{3}$''$=\omega_z$.
These last expressions only have a symbolic meaning where we interpret  
``$dq_{i}$'' and ``$\frac{\partial}{\partial q_{i}}$'', $1\leq i\leq 3$, as
adequate combinations of the differentials and partial derivatives,
respectively, of the eulerian angles. Note that
$\{\frac{\partial}{\partial x},
\frac{\partial}{\partial y}, \frac{\partial}{\partial z}$,
``$\frac{\partial}{\partial q_{1}}$'',
``$\frac{\partial}{\partial q_{2}}$'', ``$\frac{\partial}{\partial
q_{3}}$''$\}$ and $\{ dx, dy, dz,$ ``$dq^{1}$'', ``$dq^{2}$'',
``$dq^{3}$'' $\}$ are dual bases. Moreover, observe that, from the
nonintegrability of the constraints, the differential forms
``$dq^i$'' do not represent exact differentials. 

The non-holonomic generalized differentiable codistribution $D$ is
given by 
\[ 
D_{(x,y,\phi,\theta,\psi)}=\left\{
\begin{array}{cc}
\{0\} \, , & \hbox{if} \;\; x\leq 0 \, , \\
\hbox{span} \,\{ dx-rdq^2, dy+rdq^1\} \, , & \hbox{if} \; \; x>0 \, .
\end{array}
\right.
\]
The intersection of the regular set of the generalized codistribution
and the $(x,y)$-plane has two 
connected components, the half-planes $R_1=\{x<0\}$ and $R_2=\{x>0\}$. The 
line $\{x=0\}$ belongs to the singular set of $D$.

On $R_1$ the codistribution is zero, so the motion equations are
\begin{eqnarray}\label{em1}
m\ddot{x}&=&0\nonumber \, , \\
m\ddot{y}&=&0\nonumber \, , \\
mk^2\dot{\omega_x}&=&0 \, , \\
mk^2\dot{\omega_y}&=&0\nonumber \, , \\
mk^2\dot{\omega_z}&=&0\nonumber \, .
\end{eqnarray}

On $R_2$ we have to take into account the constraints to obtain the
following equations of motion
\begin{eqnarray}\label{em2}
m\ddot{x}&=&\lambda_1\nonumber \, , \\
m\ddot{y}&=&\lambda_2\nonumber \, , \\
mk^2\dot{\omega_x}&=&r\lambda_2 \, , \\
mk^2\dot{\omega_y}&=&-r\lambda_1\nonumber \, , \\
mk^2\dot{\omega_z}&=&0\nonumber \, ,
\end{eqnarray}
with the constraint equations $\dot{x}-r\omega_y=0$ and
$\dot{y}+r\omega_x=0$. One can compute the Lagrange multipliers by an
algebraic procedure described in \cite{noholo}.

Suppose that the sphere starts its motion at a point of $R_1$ with the
following initial conditions at time $t=0$:
$x_0<0$, $y_0$, $\dot{x}_0>0$, $\dot{y}_0$, $(\omega_x)_0$, 
$(\omega_y)_0$ and $(\omega_z)_0$. Integrating equations (\ref{em1}) we 
have that if $x(t)<0$
\begin{eqnarray}\label{em3}
x(t)&=&\dot{x}_0t+x_0\nonumber \, , \\
y(t)&=&\dot{y}_0t+y_0\nonumber \, , \\
\omega_x(t)&=&(\omega_x)_0 \, , \\
\omega_y(t)&=&(\omega_y)_0\nonumber \, , \\
\omega_z(t)&=&(\omega_z)_0\nonumber \, .
\end{eqnarray}

At time $\bar{t}=-x_0/\dot{x}_0$ the sphere finds the rough surface of 
the plane, where the codistribution is no longer zero and it is
suddenly forced to roll without sliding (Case 1). Following the
discussion in Section 4, we calculate the instantaneous change
of velocity (momentum) at $x=0$. 

First of all we compute the matrix ${\cal C}$
\begin{eqnarray*}
{\cal C} & = &
\left(\begin{array}{ccccc}
1&0&0&-r&0\\
0&1&r&0&0\\
\end{array}\right)
\left(
\begin{array}{ccccc}
1&0&0&0&0\\
0&1&0&0&0\\
0&0&k^{-2}&0&0\\
0&0&0&k^{-2}&0\\
0&0&0&0&k^{-2}\\
\end{array}
\right)
\left(
\begin{array}{cc}
1&0\\
0&1\\
0&r\\
-r&0\\
0&0
\end{array}
\right)
\\
\bigskip
&=&
(1+r^{2}k^{-2})
\left(
\begin{array}{cc}
1&0\\
0&1\\
\end{array}
\right) \, .
\end{eqnarray*}
Next, a direct computation shows that the projector ${\cal P}$ does not depend on the base point
\[
{\cal P} = 
\left(
\begin{array}{ccccc}
\frac{r^{2}}{r^{2}+k^{2}}&0&0&\frac{r}{r^{2}+k^{2}}&0\\
0&\frac{r^{2}}{r^{2}+k^{2}}&-\frac{r}{r^{2}+k^{2}}&0&0\\ 
0&\frac{-rk^{2}}{r^{2}+k^{2}}&\frac{k^{2}}{r^{2}+k^{2}}&0&0\\
\frac{rk^{2}}{r^{2}+k^{2}}&0&0&\frac{k^{2}}{r^{2}+k^{2}}&0\\
0&0&0&0&1
\end{array}
\right) \, .
\]
Therefore, we have
\begin{eqnarray*}
(p_{x})_{+} & = & \frac{r^{2} (p_{x})_{0} +
r (p_{2})_{0}}{r^{2}+k^{2}} \; ,\\
(p_{y})_{+} & = & \frac{r^{2} (p_{y})_{0} -
r (p_{1})_{0}}{r^{2}+k^{2}} \; ,\\
(p_{1})_{+} & = & \frac{-rk^{2} (p_{y})_{0} +
k^{2} (p_{1})_{0}}{r^{2}+k^{2}} \; ,\\
(p_{2})_{+} & = & \frac{rk^{2} (p_{x})_{0} +
k^{2} (p_{2})_{0}}{r^{2}+k^{2}} \; ,\\
(p_{3})_{+} & = & (p_{3})_{0} \; .
\end{eqnarray*}

Now, using the relation between the momenta and the quasi-velocities
\[
p_{x} = \dot{x}\, ,
p_{y} = \dot{y}\, ,
p_{1} = k^{2} \omega_{x}\, ,
p_{2} = k^{2} \omega_{y}\, ,
p_{3} = k^{2} \omega_{z} \, ,
\]
we deduce that
\begin{eqnarray}\label{em4}
\dot{x}_+&=&\frac{r^2 \dot{x}_0+rk^2(\omega_y)_0}{r^2+k^2} 
\nonumber \,
, \\
\dot{y}_+&=&\frac{r^2 \dot{y}_0-rk^2(\omega_x)_0}{r^2+k^2}\nonumber \,
, \\
(\omega_x)_+&=&\frac{-r \dot{y}_0+k^2(\omega_x)_0}{r^2+k^2} \, , \\
(\omega_y)_+&=&\frac{r \dot{x}_0+k^2(\omega_y)_0}{r^2+k^2}\nonumber \,
, \\
(\omega_z)_+&=&(\omega_{z})_{0} \nonumber \, .
\end{eqnarray}

Finally, integrating equations (\ref{em2}) at time $\bar{t}=-x_0/\dot{x}_0$ 
with initial conditions given by (\ref{em4}) we obtain that if $t>\bar{t}$
\begin{eqnarray}\label{em5}
x(t)&=&
\frac{r^2 \dot{x}_0+rk^2(\omega_y)_0}{r^2+k^2}(t-\bar{t})\nonumber \, ,
\\ 
y(t)&=&\frac{r^2 \dot{y}_0-rk^2(\omega_x)_0}{r^2+k^2}(t-\bar{t})+\dot{y}_0 
\bar{t}+y_0\nonumber \, , \\
\omega_x(t)&=&\frac{-r \dot{y}_0+k^2(\omega_x)_0}{r^2+k^2} \, , \\
\omega_y(t)&=&\frac{r \dot{x}_0+k^2(\omega_y)_0}{r^2+k^2}\nonumber \, ,
\\
\omega_z(t)&=&(\omega_z)_0\nonumber \, .
\end{eqnarray}

\subsection{Particle with constraint}

Let us consider the motion of a particle of mass $1$ in $\R^3$ subjected to the following constraint
\[
\phi = (y^2-x^2-z)\dot{x}+(z-y^2-xy)\dot{y}+x\dot{z} = 0 \, .
\]
In addition, let us assume that there is a central force system
centered at the point $(0,0,1)$ with force field given by
\[
F = -xdx-ydy+(1-z)dz \, .
\]
Then, the Lagrangian function of the particle is
\[
L = T-V = \frac{1}{2}(\dot{x}^2+\dot{y}^2+\dot{z}^2 + x^2+y^2+z^2-2z)
\, ,
\]
and the constraint defines a generalized differentiable codistribution
$D$, whose singular set is $S=\{(x,y,z):x=0,z=y^2\}$.

On $R$, the regular set of $D$, the dynamics can be computed 
following the standard symplectic procedure \cite{noholo} to obtain $\Gamma_{L,D} = \Gamma_L + \lambda Z$, where
\begin{eqnarray*}
\Gamma_L &=& \dot{x}\frac{\partial}{\partial x} + \dot{y}\frac{\partial}
{\partial y} + \dot{z}\frac{\partial}{\partial z} - x\frac{\partial}{\partial 
\dot{x}} - y\frac{\partial}{\partial \dot{y}} -
(z-1)\frac{\partial}{\partial \dot{z}} \, , \\
Z&=& -\left( (y^2-x^2-z)\frac{\partial}{\partial \dot{x}} + (z-y^2-xy) 
\frac{\partial}{\partial \dot{y}} + x\frac{\partial}{\partial \dot{z}}
\right) \, ,
\end{eqnarray*}
and $\lambda$ is given by
\begin{equation}\label{madrecordero}
\lambda=-\frac{\Gamma_L(\phi)}{Z(\phi)} =
\frac{-2x\dot{x}^2+y\dot{y}\dot{x}-2y\dot{y}^2-x\dot{y}^2+\dot{y}\dot{z}+x^3+
y^3-yz+ x}
{(y^2-x^2-z)^2+(z-y^2-xy)^2+x^2} \, .
\end{equation}
Consequently, the motion equations on $R$ are
\begin{eqnarray}\label{loco}
m\ddot{x}+x & = & \lambda(y^2-x^2-z) \nonumber \, , \\
m\ddot{y}+y & = & \lambda(z-y^2-xy) \, , \\
m\ddot{z}+z-1 & = & \lambda x \nonumber \, ,
\end{eqnarray}
with the constraint equation $\phi=0$. 

From the discussions of \cite{Chinos}, we know that in this case there
is an integral surface, $C$, of the constraint $\phi$, that is, a surface on 
which all motions satisfy the constraint. This surface is
\[
C=\{(x,y,x): z-x^2-y^2+xy=0 \} \, .
\]
Note that $S \subset C$. Therefore, if a motion takes place on the
cone-like surface $C$, it is confined to stay on this critical surface, 
unless it reaches a singular point. In this case, the space of allowable 
motions is suddenly increased (in fact, $T\R^3$), and the motion can ``escape''
from $C$. In addition, this proves that the unique way to pass from one
point of the exterior of the $C$ to the interior, or viceversa,
is through the singular set $S$.

In particular, we are interested in knowing
\begin{enumerate}
\item Is there any trajectory satisfying equations (\ref{loco}) which passes 
through the singular set?
\item if so, which are the possible momentum jumps due to the changes in
the rank of the codistribution $D$?
\end{enumerate}

So far, we do not know an answer for the question of 
the existence of a motion of (\ref{loco}) crossing $S$. It seems that
on approaching a singular point, the constraint force can become
increasingly higher (\ref{madrecordero}). Consequently, this force possibly 
``disarranges'' the approaching of the motion to $S$. Numerical
simulations are quite useless in this task, because of the special
nature of the problem: the hard restriction given by the fact that a motion 
crossing the cone-like surface $C$ must do it through the singular 
part $S$. Indeed, the numerical simulation performed in \cite{Chinos} 
crosses the surface $C$ through points which are not in $S$.

Concerning the second question, let us suppose that there is a trajectory of 
the dynamical system (\ref{loco}), $q(t)=(x(t),y(t),z(t))$, that passes 
through a singular point at time $t_0$, i.e. $x(t_0)=0$ and $z(t_0)=y^2(t_0)$. 
The rank of the codistribution $D$ at the immediately preceding and posterior points is 1, meanwhile at $q(t_0)$ it is 0 (Case 3). So, a possible jump of the momentum can be induced by the change in the rank of $D$.

A direct computation shows that the projector ${\cal P}$ depends explicitly 
on the base point $q\in Q$. Equivalently, we have that $D^+_{q(t_0)}$ depends strongly on the trajectory $q(t)$. In fact, taking two curves $q_1(t)$, $q_2(t)$ passing through $q(t_0)$ at time $t_0$, and satisfying $x_1(t) \ll z_1(t)-y_1^2(t)$ and $z_2(t)-y_2^2(t) \ll x_2(t)$ when $t\rightarrow t_0^+$ respectively, one can easily see that $D^+_{q_1(t_0)} \ne D^+_{q_2(t_0)}$ (the expression $f(t) \ll g(t)$ when $t\rightarrow t_0^+$ means that $\lim_{t\rightarrow t_0^+} f(t)/g(t) = 0$).

Consequently, we are not able to give an answer to question (ii) (in case the first one was true) unless we assume some additional information: for example, that the balance between $x(t)$ and $z(t)-y^2(t)$ is the same for 
$t\rightarrow t_0^-$ and $t\rightarrow t_0^+$. In such a case, 
$D^-_{q(t_0)}=D^+_{q(t_0)}$ and we would conclude that there is no jump. In 
mechanical phenomenae of the type sliding-rolling, as the ones studied in
Section 5.1, this kind of ``indeterminacy" will not occur in general. 

In spite of the fact that the most natural thing in this case seems to be to think that there is no jump of momenta, a mathematical explanation of it is still to be found.

\section*{Acknowledgements}

This work was partially supported by grant DGICYT (Spain) PB97-1257.
J. Cort\'es and S. Mart\'{\i}nez wish to thank Spanish Ministerio de 
Educaci\'on y Cultura for FPU and FPI grants, respectively. The authors
wish to thank F. Cantrijn for helpful comments and suggestions. We also wish to thank the referees for their valuable criticism which contributed to improve this paper.

{\parindent 0cm

{\sc Jorge Cort{\'e}s \dag \quad Manuel de Le\'on \ddag \quad Sonia Mart{\'\i}nez \S} 

{\it Instituto de Matem\'aticas y F{\'\i}sica Fundamental \\
Consejo Superior de Investigaciones Cient{\'\i}ficas\\
Serrano 123, 28006 Madrid, SPAIN\\
e-mail: \dag j.cortes@imaff.cfmac.csic.es \quad
\ddag mdeleon@imaff.cfmac.csic.es \\
\S s.martinez@imaff.cfmac.csic.es}

{\sc David Mart{\'\i}n de Diego}

{\it Departamento de Econom{\'\i}a Aplicada (Matem\'aticas)\\
Facultad de CC. Econ\'omicas y Empresariales\\
Universidad de Valladolid\\
Avda. Valle Esgueva 6, 47011 Valladolid, SPAIN\\
e-mail: dmartin@esgueva.eco.uva.es}      

\end{document}